\documentclass[11pt]{article}
\usepackage{amsfonts}
\usepackage{amsmath}
\usepackage{amssymb}
\usepackage{tikz}
\newcommand{\node}[2]{\put(#1,#2){\circle*{5}}}
\def\R{\mathbb{R}}
\def\N{\mathbb{N}}

\def\H{\mathcal{H}}
\def\nl{\newline}

\def\be{\begin{eqnarray*}}
\def\en{\end{eqnarray*}}
\def\eps{\varepsilon}

\tikzset{
  hello/.style={
    label={[draw,circle,xscale={1.25},minimum size=14mm+#1*21mm]center:},
    label={[red,yshift=#1*11mm]center:hello}
  }
}

\begin{document}
\begin{center}
{\bf {\Large What lies between polynomial and exponential growth?}
\vspace{0.2in}\\ Titus Hilberdink}
%Department of Mathematics, University of
%Reading, Whiteknights,\\
%PO Box 220, Reading RG6 6AX, UK; t.w.hilberdink@reading.ac.uk
\vspace{0.15in}\\
\end{center}

\indent
\begin{abstract} 
In this paper we give an alternative exposition of a recent paper regarding the classification of growth rates of real functions. We take a different point of view, focussing on understanding possible growth rates between polynomial and exponential. In order to be able to explicitly name a range of such functions, we first need to extend our basic functions. We do this via a `tower' of Abel functions. With these one can classify functions in a natural way with  polynomials and exponentials in consecutive classes. We show there are large gaps between these classes which indicate that it is mostly unknown what lies between polynomial and exponential growth, especially if the ``Continuum Hypothesis for classes'' is true.

%In this paper we discuss the possible rates of growth of functions between polynomial and exponential growth. In order to be able to %explicitly name a range of such functions beyond the usual log-exp functions (like $x^{\log x}, e^{\sqrt{x}}$) we need to extend our %basic functions. We do this via a `tower' of Abel functions, the first of these satisfying $H(e^x)=H(x)+1$. Armed with these, one can %explicitly write down many examples and also
%classify functions in a natural way according to their rate of growth. This raises further questions about the nature of these classes %and how many there are.   Conclusion: it is mostly unknown what lies between polynomial and exponential, especially if the %``Continuum Hypothesis'' regarding classes is true.
%In this paper we develop a classification of real functions based on growth rates of repeated iteration. We show how functions %are naturally distinguishable when considering inverses of repeated iterations. For example, $n+2\to 2n\to 2^n\to %2^{\cdot^{\cdot^2}}$ ($n$-times) etc. and their inverse functions $x-2, x/2, \log x/\log 2,$ etc. Based on this idea and some %regularity conditions we define classes of functions, with $x+2$, $2x$, $2^x$ in the first three classes. 
%
%We prove various properties of these classes which reveal their nature, including a `uniqueness' property.  We exhibit examples %of functions lying between consecutive classes and indicate
%how this implies these gaps are very `large'. Indeed, we suspect the existence of a {\em continuum} of such classes.\nl
\medskip

\vspace{0.3in}\end{abstract}

\begin{center}
{\bf {\large \S 1. Introduction}}
\end{center}
In this article, we explore what possible growth rates lie between polynomial and exponential growth. By this we mean functions $f$ defined on some interval $(a,\infty)$ which tend to infinity faster than polynomial but slower than exponential; i.e. \[ x^A\ll f(x)\ll e^{\eps x}\]
for all $A,\eps>0$. In particular, we are interested in explicitly naming functions which indicate the full spectrum of possible growth rates. 
Of course, one can readily produce examples explicitly, functions like $x^{\log x}$ and $e^{\sqrt{x}}$, but are we really any closer to bridging the gap? G. H. Hardy, in his monograph {\em Orders of Infinity} \cite{H}, introduced his class $\H$ of {\em logarithmico-exponential} functions, or log-exp functions for short, which are the functions obtained from a finite number of applications of the operations $+,-,\times,\div,\exp$ and $\log$ on the constant functions and $x$. For functions of this class he proved some remarkable properties. In particular, regarding their behaviour at infinity he proved:
%
%(i) $\H$ is closed under differentiation;
%
%(ii) every non-zero function in $\H$ is eventually either positive or negative;
%(iii) 
{\em for $f\in\H$ such that $f(x)\to\infty$ as $x\to\infty$, there exist $r,s\in\N$ and $\mu >0$ such that for all $\delta>0$,
\[(\log_s x)^{\mu-\delta}< \log_r f(x)< (\log_s x)^{\mu+\delta}\tag{1.1}\]
for $x$ sufficiently large}\footnote{Here $\log_n x=\log(\log_{n-1}x)$ and $\log_1 x=\log x$, while $e_k^x=\log_{-k}x$ denotes the exponential function iterated $k$ times. Also, we write $f\sim g$ if $f(x)/g(x)\to 1$ as $x\to\infty$ and $f\prec g$ if $f(x)/g(x)\to 0$.} (see \cite{H}). Hence $\log_{r+1} f(x)\sim \mu\log_{s+1} x$ as $x\to\infty$, so that 
\[\log_{r+2} f(x)\sim \log_{s+2} x.\]
As such, $\log_{r+n} f(x)=\log_{s+n} x+o(1)$ for each $n\ge 3$. Call $r-s$ the {\em order} of $f$.  Thus $e^x$ has order 1, $e^{e^x}$ has order 2 and $\log_k x$ has order $-k$.
So for $f\in\H$ between polynomial and exponential, $f$ has order 0 or 1 and, in this sense, these log-exp functions cannot really bridge the gap. The largest log-exp functions of order 0 are amongst the sequence of functions $e_k^{\log_k x+1}$ for $k=0,1,2,\ldots$; i.e
\[ x+1, e^{\log x +1}=ex,e^{e^{\log\log x+1}}=x^e, e^{e^{e^{\log\log\log x+1}}}=e^{(\log x)^e}, e^{e^{e^{e^{\log\log\log\log x+1}}}},\ldots\tag{1.2}\]
while the smallest log-exp functions of order 1 are amongst the sequence of functions $e_{k+1}^{\log_k x-1}$.

Hardy noticed these `gaps' between the growth rates of log-exp functions and, as an example of a function lying outside $\mathcal{H}$,  he discussed a function $\varphi$ which satisfies \[\varphi(\varphi(x))=e^x.\]
We can call such a function a $\frac12^{\rm th}$-iterate of $e^x$. Notice that if we also have\footnote{Indeed this is forced if $\varphi$ is continuous.} $x\prec \varphi(x)\prec e^x$, then (using the fact that $\log_n \varphi(x) = \varphi(\log_n x)$) we have
\[  \log_n x\prec \log_n \varphi(x)\prec \log_{n-1} x\quad\mbox{ for every $n$.}\]
This shows $\varphi$ is larger than all order 0 log-exp functions and smaller than all order 1 log-exp functions. 
More generally, this is the case for other fractional iterates of $e^x$.

It may look like this function $\varphi$ lies ``halfway'' between $x$ and $e^x$ but, in a certain sense, $\varphi(x)$ is not so far from $e^x$. From the point of view of composition, $\varphi(x)$ is to $e^x$ what $2x$ is to $4x$, or $x^{\sqrt{5}}$ is to $x^5$. These latter examples one typically lumps together, the first two being both linear, the last two of polynomial type. We shall see that, even by including functions like $\varphi$ and other positive fractional iterates of $e^x$, there still remains a large gap but this is less obvious. In order to see this and get a better view and understanding of all the possible growth rates between $x$ and $e^x$, one needs more basic functions to work with.  These we introduce in section 2.

We briefly summarize the rest of the paper. In \S2, we introduce two ideas; a sequence of basic functions obtained via the Abel functional equation and the order of a function with respect to other functions. With these concepts we construct classes of functions in \S3. We discuss some relevant properties of these classes, in particular the classes containing the polynomial and exponential functions respectively. 

In \S4 we look at functions between these classes and indicate how this leads naturally to a large gap and, heuristically, to the existence of a continuum of classes. In \S5, we tie this in with the existence of a suitable extension of the Ackermann function. 

\medskip

The results about classes in this article are essentially all proven in \cite{TH3}. However, here we take a different point of view which gives perhaps a more readily accessible explanation of the concepts, their motivation, and significance of the results. 
\nl

\medskip

\noindent
{\bf Notation}\, Unless stated otherwise, all functions are considered
to be defined on a neighbourhood of infinity. By `$f$ is continuous/increasing/etc.' 
we mean `$f$ is continuous/increasing/etc. on some interval $[A,\infty)$'. 
Also we write $f<g$ to mean $\exists\, x_0$ such that for $x\geq x_0$, $f(x)<g(x)$.

We have the usual definitions for $f\sim g$, $f=o(g)$, $f\prec g$, $f\succ g$, namely
$f(x)/g(x)$ tends to $1,0,0,\infty$ respectively, as $x\to\infty$. By $f\asymp g$, 
we mean $\exists\, a,b>0$ such that $a<f(x)/g(x)<b$ on a neighbourhood of infinity. 

For a given $f$, we write $f^n$ for the $n^{\rm th}$-iterate. For the special functions $\exp$ and $\log$ we use the notation $e_n^x$ and $\log_n x$ for the $n^{\rm th}$-iterates at $x$ respectively. 

Also let ${\rm D}_{\infty}^+ $ denote the class of $C^1$ functions defined on some neighbourhood of infinity tending to infinity and whose derivative is positive;  i.e.
\[ {\rm D}_{\infty}^+ =\{f:[A,\infty)\to\R:f\mbox{ is continuously differentiable, $f^{\prime}>0$ and $f(x)\to\infty$ as $x\to\infty\}$.}\]
%\end{align*}
Note that ${\rm D}_{\infty}^+$ is a group under composition, if we identify functions which are equal on a neighbourhood of infinity.\nl
\bigskip

\begin{center}
{\bf {\large \S 2. The Abel equation and the functions $\Xi_n(x)$}}
\end{center}
The function $\varphi(x)$ and other iterates of $e^x$ mentioned above can be obtained directly as follows: let $H$ be a strictly increasing and continuous function satisfying the {\em Abel functional equation} for $e^x$:
\[ H(e^x)=H(x)+1.\tag{2.1}\]
Such a function is easily shown to exist. Indeed, more generally, given a function $f\in {\rm D}_{\infty}^+$ such that $f(x)>x$, there are always ${\rm D}_{\infty}^+$ solutions of 
\[F(f(x))=F(x)+1.\]
We shall call such $F$ an {\em Abel function of $f$}. For example, to obtain a strictly increasing and continuous solution, we just need $F$ continuous and strictly increasing on $[a, f(a)]$ (with $a$ suitably large) such that $F(f(a))=F(a)+1$. Then extend $F$ to $[a,\infty)$ via $F(f(x))=F(x)+1$. For the existence of ${\rm D}_{\infty}^+$ solutions see the appendix.
For more information on the Abel equation, see for example \cite{KCG}.

Now with $H$ as in (2.1), the function $\varphi(x)=H^{-1}(H(x)+\frac12)$ satisfies
%\footnote{Indeed, every such $\varphi$ is of this form.} 
$\varphi(\varphi(x))=e^x$. More generally, $H^{-1}(H(x)+\frac{p}{q})$ with $p,q\in\N$ is a $(\frac{p}{q})^{\rm th}$-iterate of $e^x$. Notice that $H(x)$ tends to infinity very slowly --- slower than any iterate of $\log$; equivalently, $H^{-1}(x)$ tends to infinity faster than any iterate of $e^x$. 

\medskip

\noindent
{\bf 2.1 The functions} $\Xi_n(x)$. We shall need to name more functions. Fixing one ${\rm D}_{\infty}^+$ solution $H$ of the Abel equation for $e^x$, we can now do the same for $H^{-1}$; i.e. obtain a ${\rm D}_{\infty}^+$ Abel function of $H^{-1}$ (note that $H^{-1}(x)>x$ for large enough $x$). Then we can repeat this indefinitely. In this way, we obtain a sequence of ${\rm D}_{\infty}^+$ functions $\Xi_n$ satisfying 
\[ \Xi_n(x)=\Xi_n(\Xi_{n-1}(x))+1\qquad\mbox{ for $n\in\N_0$}.\]
That is, $\Xi_n$ is an Abel function of $\Xi_{n-1}^{-1}$. 
In fact we will start further down the list and take $\Xi_0(x)=x-e$, $\Xi_1(x)=\frac{x}{e}$, $\Xi_2(x)=\log x$ and 
$\Xi_3(x)=H(x)$ as above. Observe that, as for the case $n=3$, $\Xi_n(x)\to\infty$ as $x\to\infty$ slower than any iterate of $\Xi_{n-1}(x)$. Actually, we can choose the $\Xi_n$ so that they are concave but we will not make use of this here.  
 As in \cite{TH3}, we say a function $f$ is of {\em finite class} if there exists $k$ such that $\Xi_k<f<\Xi_k^{-1}$. 

\medskip

In a sense, the derivatives of the $\Xi_n$ are more important. Note that $\Xi_0^\prime(x)=1$, $\Xi_1^\prime(x)=\frac{1}{e}$, $\Xi_2^\prime(x)=\frac1{x}$. Writing $\Xi_n^\prime(x)=\frac1{\chi_n(x)}$, we have 
\[ \chi_n(x) = \chi_{n-1}(x)\chi_n(\Xi_{n-1}(x)).\tag{2.2}\]
For example, $\chi_3(x) = x\chi_3(\log x)$. Thus, for large $x$, $\chi_3(x)\asymp x\log x\log\log x \cdots \log_k x$ where $k$ is such that $\log_k x\in [1,e]$. The functions $\chi_n$ with $n\ge 4$ are marginally larger.  In any case, they are all $\ll x(\log x)^2$. This behaviour of $\chi_n(x)$ for large $x$ should not come as a surprise when we realize that 
\[ \int_a^\infty \frac1{\chi_n(x)}\, dx=\infty\quad\mbox{ while }\quad \int_a^\infty \frac1{\chi_n(x)\Xi_n(x)^2}\, dx<\infty.\]
With these new functions, we can explicitly write down many other functions lying between polynomial and exponential growth. 
Even just using the function $H$ (i.e. $\Xi_3$), we can write down many. First note that with $f_k(x) =e_k^{\log_k x+1}$ as in (1.2), we have (for some $\theta_x\in (0,1)$)
\[ \Xi_3(f_k(x)) = \Xi_3(\log_k x+1)+k = \Xi_3(\log_k x)+k+\frac{1}{\chi_3(\log_k x+\theta_x)}= \Xi_3(x)+\frac{1}{\chi_3(\log_k x+\theta_x)}.\] 
As these $f_k$ are the largest functions in $\mathcal{H}$ of order 0, this shows that $\Xi_3(f(x))-\Xi_3(x)$ must go to zero at least as fast as $\frac1{\log_k x}$ for some $k$ for every log-exp function $f$ of order 0. By comparison, for $\varphi$, $\Xi_3(\varphi(x)) = \Xi_3(x)+\frac{1}{2}$, while the ${\frac1n}^{\rm th}$-iterate of $e^x$, say $\varphi_n$, satisfies $\Xi_3(\varphi_n(x)) = \Xi_3(x)+\frac{1}{n}$. 

\medskip

For example, consider the functions
\[ g(x)=\Xi_3^{-1}\Bigl(\Xi_3(x)+\frac1{\Xi_3(x)}\Bigr) \qquad\mbox{ and }\qquad h(x)=\Xi_3^{-1}\Bigl(\Xi_3(x)+\frac1{\Xi_3^{-1}(\frac12\Xi_3(x))}\Bigr).\tag{2.3}\]
Notice that for $g$, $\Xi_3(g(x))-\Xi_3(x)$ tends to $0$ like $\frac1{\Xi_3(x)}$ so it is larger than any order 0 log-exp function, yet it is smaller than any positive fractional iterate of $e^x$. Similarly for $h$ because $\Xi_3(h(x))-\Xi_3(x)$ tends to zero and $\Xi_3^{-1}(\frac12\Xi_3(x))<\log_k x$ for any $k$ for $x$ large enough. Thus both functions lie between all the $f_k$ and the $\varphi_n$. 

For an example using $\Xi_4$, consider
\[ \ell(x) = \Xi_3^{-1}\Bigl(\Xi_3(x)+\frac1{\Xi_4^{-1}(\Xi_4(x)-\frac12)}\Bigr).\tag{2.4}\]
Note that $\Xi_4^{-1}(\Xi_4(x)-\frac12)$ is a $\frac12^{\rm th}$-iterate of $\Xi_3(x)$. One can check that $h<\ell<g$.
We will say more about these examples in section 3. 

\bigskip

\noindent
{\bf 2.2 Orders of functions.} \, Another concept we require is the notion of order of a function with respect to another function (see \cite{TH1}). Given a strictly increasing continuous and unbounded function $F$, we define $O_F(f)$ --- the {\em order of $f$ with respect to $F$} --- by
\[ O_F(f)=\lambda\qquad \mbox{ if }\qquad \lim_{x\to\infty} \Bigl(F(f(x))-F(x)\Bigr) =\lambda\qquad (\lambda\in\R).\]
For example, with $F(x)=\log x$, $O_F(f)=\lambda\iff f(x)\sim e^\lambda x$; with $F(x)=\log\log x$, $O_F(f)=\lambda\iff f(x) = x^{e^\lambda+o(1)}$. So, crudely speaking, orders w.r.t. $\log$ distinguish between multiples of $x$ while those for $\log\log$ distinguish between powers of $x$.

Notice that if $f\in\mathcal{H}$ of order $k$ (as defined in \S 1), then $\Xi_3(f(x))-\Xi_3(x)\to k$ since, for $n$ large enough,
\[ \Xi_3(f(x))-\Xi_3(x) = \Xi_3(\log_nf(x))-\Xi_3(\log_n x) =  \Xi_3(\log_{n-k}x+o(1))-\Xi_3(\log_n x)=k+o(1).\]
Thus this notion generalizes the earlier notion of order. Also the functions in (2.3) and (2.4) have order 0 w.r.t. $\Xi_3$ while $\varphi$ has order $\frac12$. 
\medskip

Observe that with $\lambda=1$, $O_F(f)=1$ says $F(f(x))=F(x)+1+o(1)$, so $F$ is an approximation to an Abel function of $f$. Indeed $F\sim G$ for any Abel function $G$ of $f$.

Finally, we point out that the general rule $O_F(f\circ g) = O_F(f)+O_F(g)$. For a detailed discussion about such orders and their use in {\em uniqueness} of fractional iterates, see \cite{TH1}.

\bigskip

\begin{center}
{\bf {\large \S 3. Classes $C_n$}}
\end{center}
Now that we have the functions $\Xi_n$, we can generate many more functions about which we can, with some justification, say we understand its growth rate. Each $\Xi_k$ can be explicitly calculated and its long term behaviour is determined by $\Xi_{k-1}$ via $\Xi_k(\Xi_{k-1}^{-n}(a)) = n+\Xi_k(a)$.

In \cite{TH3}, based on the sequence $\Xi_n$, we formed classes of functions of similar growth to each $\Xi_n^{-1}$. To describe the procedure, we first restrict to functions whose derivatives behave regularly w.r.t. this sequence. Let 
\[ \mathcal{B} = \{ f\in {\rm D}_{\infty}^+: (\Xi_n\circ f)^\prime \sim \Xi_n^\prime \mbox{ for some $n$}.\}\]
Notice that for a function $f$ of finite class there must be some $n$ such that\footnote{For $\Xi_k<f<\Xi_k^{-1}$ implies $\Xi_{k+1}-1<\Xi_{k+1}(f)<\Xi_{k+1}+1$, so $\Xi_{k+1}(f)\sim \Xi_{k+1}$.} $\Xi_n(f(x))\sim \Xi_n(x)$. Now $f\in\mathcal{B}$ just means that this asymptotic equivalence can be differentiated: $(\Xi_n(f(x)))^\prime\sim \Xi_n^\prime(x)$.

Observe that $\Xi_n\in\mathcal{B}$ for all $n$ and $\mathcal{B}$ is a group under composition. 

\medskip

We build up the aforementioned classes recursively. We start with just one function in each class, namely $\Xi_n^{-1}$; i.e. let 
\[ C_n^{(0)} = \{\Xi_n^{-1}\}\qquad (n\in\N_0).\]
For the next step, we add some functions `close' to $\Xi_n^{-1}$ in the following way: we take all the $f\in\mathcal{B}$ such that $O_{\Xi_{n+1}}(f)=1$. In other words,
\[ C_n^{(1)} = \{ f\in\mathcal{B}: O_{\Xi_{n+1}}(f)=1 \}.\]
For example for $n=1$, we have $C_1^{(0)}=\{ex\}$ and $f\in C_1^{(1)}$ if $f(x)\sim ex$ and $f^\prime(x)\sim e$.  

Now, once $C_n^{(k)}$ has been defined for some $k\ge 0$ and all $n\ge 0$, we define 
\[ C_n^{(k+1)} = \{ f\in\mathcal{B}: O_F(f)=1 \mbox{ for some $F$ with $F^{-1}\in C_{n+1}^{(k)}$}\}.\]
As such $C_n^{(0)}\subset C_n^{(1)}\subset C_n^{(2)}\cdots$ and define 
\[C_n=\bigcup_{k=0}^\infty C_n^{(k)}.\]
Thus $f_0\in C_n$ means $f_0\in C_n^{(k)}$ for some $k$, which means there exist $f_1,\ldots,f_k$ with $f_r\in C_{n+r}^{(k-r)}$ such that with $F_r=f_r^{-1}$,
\[ F_r(f_{r-1}(x)) = F_r(x)+1+o(1)\qquad\mbox{ for $r=1,\ldots,k$.}\]
In other words, after $k$ steps, we are back to some $\Xi_m$ which we `know'. 
For example, $f_0(x)=x+\frac{x}{\log x}$ lies in $C_1^{(2)}$. For with $F_1(x) = \frac12(\log x)^2$, we have $O_{F_1}(f_0)=1$. Next $f_1(x)=F_1^{-1}(x) =e^{\sqrt{2x}}$ and, with $F_2=\Xi_3$, $O_{F_2}(f_1)=1$. As $\Xi_3^{-1}\in C_3^{(0)}$, we have $f_1\in C_2^{(1)}$ and $f_0\in C_1^{(2)}$.  

A few more examples are contained in the following table.\footnote{Here $a$ is a constant greater than 1. Also, by $f\simeq g+h$ we mean $f-g\sim h$.}
\[  \begin{array}{|c|cc@{\quad\vdots\quad}cccc@{\quad\vdots\quad}cc|}
\hline
f_0 & x+a  & x+\sqrt{x} & x+\frac{x}{\log x} & ax & x^a & x^{\log x}& \varphi(x) & e^x \\
\hline
F_0 & x-a  & \simeq x-\sqrt{x} & \simeq x-\frac{x}{\log x} & \frac{x}{a} & x^{1/a} & e^{\sqrt{\log x}} & \varphi^{-1}(x) & \log x  \\
F_1 & \frac{x}{a} & 2\sqrt{x} & \frac{1}{2}(\log x)^2 & \frac{\log x}{\log a} & \frac{\log\log x}{\log a} & \frac{\log\log\log x}{\log 2} & 2\Xi_3(x) & \Xi_3(x)  \\
F_2 & \frac{\log x}{\log a} & \frac{\log\log x}{\log 2} & \Xi_3(x) & \Xi_3(x) & \frac{1}{2}\Xi_3(x) & \frac13\Xi_3(x) & \Xi_4(x) & \Xi_4(x)  \\
F_3 & \Xi_3(x) & \frac{1}{2}\Xi_3(x) & \Xi_4(x) & \Xi_4(x) & \Xi_4(x) & \Xi_4(x) & \Xi_5(x)& \Xi_5(x) \\
F_4 & \Xi_4(x) & \Xi_4(x) & \Xi_5(x) & \Xi_5(x) & \Xi_5(x) & \Xi_5(x)& \Xi_6(x) & \Xi_6(x) \\
\hline
\end{array} 
 \]
From this table we can read off that the first two functions $f_0$ lie in $C_0$ (indeed they are in $C_0^{(3)}$ and $C_0^{(4)}$ respectively), the next four lie in $C_1$, and the final two lie in $C_2$. 

\medskip

Many properties of these classes were proved in \cite{TH3}, of which we mention those most relevant.
\begin{enumerate}
\item Every log-exp function greater than $x+1$ lies in $C_0^{(4)}\cup C_1^{(3)}\cup C_2^{(2)}$. 
\item Functions in $C_n$ are eventually smaller than those in $C_{n+1}$; more compactly, we write this as 
\[C_n<C_{n+1}.\]
\item For every $k\ge 1$, there are functions in $C_n^{(k+1)}$ which are larger/smaller than any function in $C_n^{(k)}$.
\end{enumerate}
These were proved in Proposition A.4, Corollary 2.4(b) and Theorem 2.7  of \cite{TH3} respectively. For the proofs of (b) and (c), the regularity of the functions (namely that they lie in $\mathcal{B}$) was essential.

\medskip

The following picture may give a rough idea of these classes.

\bigskip
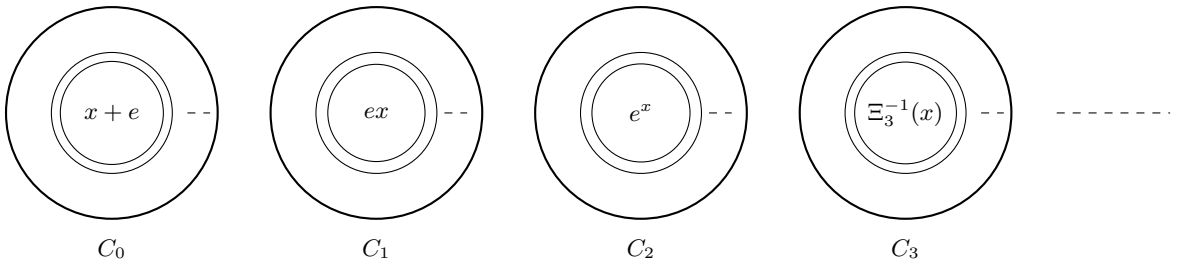
\footnotesize{
%\tikz \draw node[circle, draw]{$x+e$} ;\quad
%\tikz \draw node[circle, draw]{$ex$};\quad
%\tikz \draw (3,-3)  node[circle, draw]{$e^x$};\quad
%\tikz \draw (3,-3)  node[circle, draw]{$\Xi_3^{-1}(x)$};\qquad $\cdots$ \qquad
%\tikz \draw (3,-3)  node[circle, draw]{$\Xi_n^{-1}(x)$};$\cdots$

 \begin{tikzpicture}

 \draw (0,0) node[circle, draw]{$\phantom{a}x+e\phantom{a}$};
 \draw (0,0) circle (0.8cm);
 \draw[thick, label=above:Top] (0,0) circle (1.4cm);
\node at ( 0, -1.8) {$C_0$};
\draw (3.5,0) node[circle, draw]{$\phantom{a}\phantom{a}ex\phantom{a}\phantom{a}$};
 \draw (3.5,0) circle (0.8cm);
 \draw[thick] (3.5,0) circle (1.4cm);
\node at ( 3.5, -1.8) {$C_1$};
\draw (7,0) node[circle, draw]{$\phantom{a}\phantom{a}e^x\phantom{a}\phantom{a}$};
 \draw (7,0) circle (0.8cm);
 \draw[thick] (7,0) circle (1.4cm);
\node at ( 7, -1.8) {$C_2$};
\draw (10.5,0) node[circle, draw]{$\Xi_3^{-1}(x)$};
 \draw (10.5,0) circle (0.8cm);
 \draw[thick] (10.5,0) circle (1.4cm); 
\node at ( 10.5, -1.8) {$C_3$};
 \draw[dashed] (12.5,0)-- (14,0);
\draw[dashed] (1,0)-- (1.3,0);
\draw[dashed] (4.4,0)-- (4.7,0);
\draw[dashed] (7.9,0)-- (8.2,0);
\draw[dashed] (11.5,0)-- (11.8,0);
%\draw (14,0) node[circle, draw]{$\Xi_n^{-1}(x)$};
% \draw (14,0) circle (1cm);
% \draw[thick] (14,0) circle (1.6cm);
 \end{tikzpicture}

%\begin{tikzpicture}
%  \coordinate[hello/.list={0,1,2,3}];
%\end{tikzpicture}

%\begin{tikzpicture}
%\coordinate (O) at (0,0);
%\draw (O) circle (2.3);
%\draw (O)  circle (1.1);
%\draw (O) circle (0.35);
%\end{tikzpicture}
}
\normalsize

\bigskip

For the functions in (2.3) one can check that $g\in C_2^{(2)}$ and $h\in C_1^{(3)}$, while the function $\ell$ of (2.4) lies in $C_2^{(3)}$ and is smaller than anything in $C_2^{(2)}$. These indicate that we haven't gone very deep into the classes.

\medskip 

\begin{center}
{\bf {\large \S 4. Functions between classes}}
\end{center}
One can continue taking more complicated examples but, as we shall see, even all these functions barely make any inroads into the problem of bridging the gap between polynomial and exponential growth! For in \cite{TH3} (Theorem 2.8), it was shown that there are functions {\em between} $C_1$ and $C_2$; i.e. there exists $f$ such that {\em for all $g\in C_1$ and all} $h\in C_2$,
\[ g<f<h.\]
(Indeed, this is true between any two neighbouring classes.)

Functions with growth rate in between $C_1$ and $C_2$ have a rather alien growth rate. It is difficult to say anything more beyond the fact that such functions lie between $C_1$ and $C_2$. The question is now: {\em how much `space' is there between $C_1$ and} $C_2$ (or indeed between any two consecutive classes)? 

\medskip

Now, if $C_{n-1}<f<C_{n}$ and $F$ is strictly increasing and continuous and such that $O_F(f)=1$, then $C_n<F^{-1}<C_{n+1}$ (Proposition 2.9 of \cite{TH3}). Thus we obtain a sequence of functions $f_n$ such that $C_n<f_n<C_{n+1}$  and $O_{F_{n+1}}(f_n)=1$ for each $n\ge 0$, where $F_n=f_n^{-1}$. If ---  and so far this is hypothetical --- we can further show that we can take each $f_n\in\mathcal{B}$, then we can readily generate classes of functions {\em in between} the $C_n$ in the same way that the $C_n$ were defined: let $D_n^{(0)}=\{f_n\}$ for $n\ge 0$ and given $D_n^{(k)}$ for all $n$ and some $k\ge 0$, let
\[ D_n^{(k+1)} = \{ f\in\mathcal{B}: O_F(f)=1 \mbox{ for some $F$ with $F^{-1}\in D_{n+1}^{(k)}$}\}\qquad\mbox{ and }\quad D_n=\bigcup_{k=0}^\infty D_n^{(k)}\]
--- exactly as for the $C_n^{(k+1)}$. The same technique for showing $C_n<C_{n+1}$ now allows us to show $C_n<D_n<C_{n+1}$. So not only do we get functions between $C_1$ and $C_2$, but a whole class of functions of the same type as $C_1$ or $C_2$. 

However, it doesn't stop here. To aid our intuition, let us rename $D_n$ as $C_{n+\frac12}$. The same technique for obtaining functions between neighbouring classes can be used to find functions between $C_n$, $C_{n+\frac12}$ and $C_{n+1}$. Furthermore, again assuming we can prove they can be in $\mathcal{B}$, we can repeat the above to obtain classes between them. Let's call them $C_{n+\frac14}$ and $C_{n+\frac34}$. Repeating this process indicates that there are classes $C_{n+q}$ for any rational $q$ of the form $r/2^k$ with $r=0,1,\ldots, 2^k$ such that 
\[ C_{n+q}<C_{n+q^{\prime}} \]
for all such rationals $q,q^{\prime}$ with $q<q^{\prime}$. This shows a staggering amount of `space' between $C_1$ and $C_2$.
A more accurate picture is something like this:
\[ \stackrel{\bullet}{C_0} \qquad\stackrel{\bullet}{C_1} \qquad\stackrel{\cdot}{\phantom{A}}\qquad\stackrel{\cdot}{\phantom{A}}\qquad\stackrel{\cdot}{\phantom{A}}\qquad \stackrel{\bullet}{C_n} \qquad \stackrel{\cdot}{\phantom{A}}\]

\bigskip

\noindent
{\bf 4.1 A continuum of Classes?}\, However, this is not all. For any $\lambda>0$, one can find strictly monotonic sequences 
%$p_m$, $q_m$ 
of rationals of the form $r2^{-k}$ converging to $\lambda$.
%such that $p_m$ is strictly increasing while $q_m$ is strictly decreasing. 
By taking limits of sequences of functions in the corresponding classes, 
%$C_{n+p_m}$ and $C_{n+q_m}$,% 
we can find functions in between each of the $C_{n+q}$s. With the same assumption about such functions being in $\mathcal{B}$, we can obtain classes of functions here. In other words, these considerations suggest that there is in fact a {\em continuum} of classes $C_{n+\alpha}$ with $\alpha\in (0,1)$ in between each pair of neighbouring classes $C_n$ and $C_{n+1}$ and that $C_{n+\alpha}<C_{n+\beta}$ for $0\le\alpha<\beta\le 1$.\nl

\noindent
{\bf Continuum Hypothesis for Classes}\nl
For each $\lambda\ge 0$ there exists a class of functions $C_\lambda\subset\mathcal{B}$ satisfying:

(i) $C_\lambda<C_\mu$ whenever $\lambda<\mu$ and

(ii) $f\in C_\lambda$ implies the existence of $F$ with $F^{-1}\in C_{\lambda+1}$ such that $O_F(f)=1$.\nl

If true, this indicates that the $C_n$ $(n\in\N_0)$ represent only an {\em infinitesimally small proportion} of the functions of finite class, and that there is an astonishingly large gap of growth rates between consecutive classes --- and in particular, between polynomial and exponential. Furthermore, it is difficult to say anything meaningful about comparing the growth rates of functions from, say, $C_{\frac14}$ and $C_{\frac12}$, beyond the fact that the former are smaller than the latter.  To say something meaningful, one needs an extension of the sequence $\Xi_n$ $(n\in\N_0)$ to a continuum $\Xi_\lambda$ $(\lambda\in [0,\infty))$.\nl

There is a further interesting question of completeness: {\em can a continuum of classes be defined as above such that if $C_{\alpha-\eps}<f<C_{\alpha+\eps}$ for all $\eps>0$, then $g<f<h$ for some $g,h\in C_\alpha$?}

\bigskip

\noindent
{\bf 4.2 Big gaps between Classes.}\,  All these questions are still open, even whether $C_{n+\frac12}$ exists. However, we think this is essentially only a technical issue, namely that these sequences $f_n$ in between $C_n$ and $C_{n+1}$ lie in $\mathcal{B}$. In recent progress (unpublished), we know that we can make any {\em finite} sequence of such $f_n$ lie in $\mathcal{B}$.

  In any case, that there really is a lot of room between consecutive classes was shown in \cite{TH3} at the end of section 2.
With $f_n$ such that $C_n<f_n<C_{n+1}$ and $O_{F_{n+1}}(f_n)=1$ for each $n\ge 0$, where $F_n=f_n^{-1}$ as found above, define sets $E_n^{(k)}$ recursively as follows: first let
\[ E_n^{(0)} = \Bigl\{ f: f^{-1}=\Xi_m^{-1}\Bigl(\Xi_m-\frac{1+\delta}{\chi_m(F_{n+1})}\Bigr)\Bigr\},  \]
where $\delta$ is any continuous function decreasing to 0 and $m$ large enough, say $m\ge n+3$. As such, $C_n< E_n^{(0)} < C_{n+1}$ by Proposition 2.10 of \cite{TH3}. Now, suppose $E_n^{(k)}$ has been defined for some $k\ge 0$ and all $n\ge 1$, let
\[ E_n^{(k+1)} = \Bigl\{ f: f^{-1} = \Xi_m^{-1}\Bigl(\Xi_m - \frac{1+\delta}{\chi_m(F)}\Bigr) \mbox{ where $\delta\searrow 0$ and } F^{-1}\in E_{n+1}^{(k)} \Bigl\}. \]
With the methods in \cite{TH3}, one finds by induction that $C_n< E_n^{(k)} < C_{n+1}$ for all $k\ge 0$. Letting $E_n = \cup_{k\ge 0} E_n^{(k)}$, gives $C_n< E_n < C_{n+1}$. This shows that, at the very least, something akin to a class of functions lies between two neighbouring classes.

%As mentioned in \cite{TH3}, we could continue showing the existence of functions $g_m\in C_n^{(m)}$, $h_m\in E_n^{(m)}$ such that $g_m>C_n^{(m-1)}$ and $h_m<E_n^{(m-1)}$ and %$g_m<g_{m+1}<h_{m+1}<h_m$ to obtain functions between the $C_n, E_n$ and $C_{n+1}$. 
%
%CAN WE DO MORE HERE?

\bigskip

\begin{center}
{\bf {\large \S 5. A two-variable extension of the Ackermann function}}
\end{center}
If the continuum hypothesis for classes is correct, there would exist a continuum of functions $\{\Xi_\alpha\}_{\alpha\in [0,\infty)}$ such that
\[ \Xi_\alpha(x)=\Xi_\alpha(\Xi_{\alpha-1}(x))+1\qquad\mbox{ for $\alpha\ge 1$}\tag{5.1}\]
and further that the classes generated from $\{\Xi_{n+\alpha}\}_{n\ge 0}$ are separate. 
\medskip

The first part of this is essentially a two-variable extension of the {\em Ackermann} function (see \cite{Ack}).
% For, in terms of inverses (5.1) says 
%\[ \Xi_\alpha^{-1}(x+1) = \Xi_{\alpha-1}^{-1}(\Xi_\alpha^{-1}(x)).\tag{5.2}\]
%which one recognizes as the defining recursion for the Ackermann function. 
There are several variants of this function defined as $A:\N^2_0\to \N$ via the recursion
\[ A(m, n+1) = A(m-1,A(m,n))\qquad (m\ge 1,n\geq 0) \]
and some initial values which we will take to be $A(m,0)=2$, $A(0,n)=n+2$. As such, one has $A(1,n)=2n+2$, $A(2,n)=2^{n+2}-2$ and  $A(3,n)= 2^{\cdot^{\cdot^2}}\mbox{ ($n+2$-times)}-2$, after which one runs out of standard notation.  (This function plays an important role in recursion and computability, see for example \cite{Grz}, \cite{Rit}.)

For each fixed $m$, we can define an extension of $A(m,\cdot)$, which we will denote by $A_m:[0,\infty)\to \R$, in such a way that each such function is strictly increasing and continuous and its inverse $G_m(x) := A_m^{-1}(x)$ satifies the Abel functional equation
\[ G_m(x) = G_m(G_{m-1}(x)) +1.\tag{5.2} \]
Note the identical form in (5.1). For example, we can take $G_0(x)=x-2$, $G_1(x)=\frac{1}{2}x-1$, and 
\[ G_2(x) = \frac{\log (x+2)}{\log 2} - 2. \]
Clearly we can vary the initial values of the Ackermann function and don't need to insist that $G_m(2)=0$ or that $G_m^{-1}(n)$ is an integer. All we really want is that each $G_m\in {\rm D}_\infty^+$ and $G_m(x)\le x-a$ for some $a>1$. In other words, we could start with $G_0(x)=\Xi_0(x)=x-e$ so that we can choose $G_m=\Xi_m$ for each $m$ to satisfy (5.2). 

Now we seek a suitable extension in the second variable $m$ to real values. Suppose we have extended $\Xi_m(\cdot)$ to $\Xi_\alpha(\cdot)$ for any real number $\alpha \ge 0$. Then for each $\alpha\in [0,1)$, $\{\Xi_{n+\alpha}(\cdot)\}_{n\geq 0}$ generates classes, which we may denote by $C_{n+\alpha}$.  Now we want to make sure that $C_{\lambda}<C_{\mu}$ whenever $\lambda<\mu$. For this reason we must ensure that any such generalization satisfies the conditions
\[ (i)\quad \Xi_\alpha\in\mathcal{B}\quad\mbox{ for all $\alpha\ge 0$; }\quad (ii)\quad \Xi_{\beta}(x) = o(\Xi_{\alpha}(x))\quad\mbox{ for $\beta>\alpha\ge 1$}.\tag{5.3}\]
%EXPLAIN...Here (ii) is essential to have $C_\lambda<C_\mu$ while the regularity condition in (i) is needed for the proof (or rather, it is sufficient). 
The question is whether this is possible, and moreover, whether it can be done in such a way that the completeness we discussed in section 4.1  holds for these classes.

Notice that a naive choice like, say, 
\[\Xi_\alpha(x) = \frac{2-\alpha}{e}x^{\frac{2}{\alpha}-1}+(\alpha-1)\log x \qquad (1\le\alpha\le 2)\]
will  not do. For (ii) will be satisfied for $\alpha,\beta\in[1,2]$ but it forces $\Xi_{\alpha+1}(x)\sim c_\alpha\log\log x$ for $1<\alpha<2$ (for some $c_\alpha >0$) and (ii) fails for $\alpha,\beta\in (2,3)$ and in intervals further to the right.  Of course, we need $C_1<\Xi_\alpha^{-1}<C_2$ for each $\alpha\in (1,2)$, which rules out any such naive choice. 

\medskip

Also, there is no uniqueness to the $\Xi_\lambda$ and the $C_\lambda$ for $\lambda$ non-integral in the following sense: for any increasing bijection $\psi:[0,1)\to [0,1)$, $C_{n+\psi(\alpha)}$ ($\alpha\in [0,1)$, $n\in\N_0$) runs through the same classes but at a different rate. Thus we cannot judge the rate of growth of functions in $C_{n+\frac12}$ any better than those of, say, $C_{n+\frac23}$.
\nl

%\begin{center}
%{\bf {\large \S 6. Some new results?}}
%\end{center}
%The existence of $f\in\mathcal{B}^+$ between $C_n$ and $C_{n+1}$. \nl
%\newpage
\bigskip

\noindent
{\bf Final remarks} \nl
To add some further context for the present work, we discuss some related topics and results. 
\begin{enumerate}
\item In \cite{ERT} the authors are interested in {\em scale of functions} with which one can compare other functions. Briefly, these are subsets $\mathcal{S}$ of the class $\mathcal{C}$ of continuous functions such that $1,x\in\mathcal{S}$, and for $f,g\in\mathcal{S}$,  $f/g\to\lambda\in [0,\infty]$ (i.e. any two functions are {\em comparable}), and $(f)^\alpha (g)^\beta\in\mathcal{S}$ for all $\alpha,\beta\in\R$. They show the existence of a scale of functions $\mathcal{S}$ which is {\em irreducible},
% -- i.e. for all $f,g\in\mathcal{S}$, either $f\prec g$, $f=g$ or $f\succ g$ --- 
{\em maximal}, 
 %i.e. if $f\in\mathcal{C}$ is comparable to every $h\in\mathcal{S}$ then $f/g\to\lambda\in (0,\infty)$ for some $g\in\mathcal{S}$.
and, on the assumption of the Continuum Hypothesis, {\em dense}.
% --- for all $g\in\mathcal{C}$ with $g\succ 1$ there exist $f_1,f_2\in\mathcal{S}$ such that $1\prec f_1\prec g\prec f_2$ --- and for which any two elements of $\mathcal{S}$ are %monotonically comparable ($f/g$ is monotonic).  Also the scale $\mathcal{S}$ must be uncountable (indeed has cardinality that of the continuum). 
These are interesting existence results but they say little about the questions we are concerned with. For a start, we are discussing only functions of finite class, while in  \cite{ERT} they discuss arbitrarily large (and small) functions. Further, we are interested in naming explicit functions between $C_1$ and $C_2$, not merely the existence of a such functions.

\item There are of course many possible functions ${\rm D}_\infty^+$ functions $H$ satisfying $H(e^x)=H(x)+1$, and more generally for solutions of other Abel equations\footnote{If $H_0$ is one such function then $H(x)=H_0(x)+p(H_0(x))$ is a solution for every period one function $p$ with $p^\prime>-1$.}. This leads to many choices for each $\Xi_n$ in turn, so we end up with a different $\mathcal{B}$ and different classes. However, the growth rates will be similar. As mentioned in the appendix of \cite{TH3}, one can use a weaker regularity condition, namely $(\Xi_n \circ f)^\prime\asymp \Xi_n^\prime$ in the definition of $\mathcal{B}$ (i.e. $\asymp$ in place of $\sim$), whereby all the different Abel solutions are lumped together in the same class.

It is an interesting question if there is a particular $H$ satisfying (2.1) which is in some sense `best'. This question was explored recently in \cite{TH2}. See also the papers by Szekeres, \cite{S3}, \cite{S4}.

%\item Of course, interest in a certain type of growth rate depends on the particular topic one is interested in. For the theory of entire functions, the relevant growth rates of interest are %$e^{|z|^\rho}$ with $\rho\ge 0$. We have focussed on the gap between polynomial and exponential (or $C_1$ and $C_2$) but we could equally well have considered the gap between %$C_0$ and $C_1$, or any other neighbouring classes.

\item We also should mention that Hardy's class $\mathcal{H}$ has been extended by Boshernitzan and Rosenlicht to include solutions of a class of algebraic differential equation in \cite{Bosh}, \cite{Bosh2}, \cite{R2}. The functions obtained have similar growth properties to those in $\mathcal{H}$, namely they satisfy (1.1). These are examples of more general Hardy fields (see \cite{Shack1}).

\end{enumerate}

\bigskip

\noindent
{\bf Appendix}\nl
As we saw in the beginning of section 2, given a strictly increasing and continuous $f$ such that $f(x)>x$, we can always find a strictly increasing continuous Abel function $F$ of $f$. 
If $f\in {\rm D}_{\infty}^+$ we can make sure $F\in {\rm D}_{\infty}^+$ too as follows.

Let $a$ be large enough so that on $[a,\infty)$, $f$ is continuously differentiable, $f^\prime>0$ and $f(x)>x$. Define $G$ on $[a,f(a)]$ such that $G$ is positive and continuous and $G(f(a))f^\prime(a)=G(a)$. Then extend $G$ to $[a,\infty)$ via 
\[ G(f(x))f^\prime(x)=G(x) \quad (x\ge a).\]
As such, $G$ is positive and continuous here. Now let $F(x) = c\int_a^x G$ for $x\ge a$ with $c$ chosen so that $F(f(a))=1$ (i.e. $1/c = \int_a^{f(a)}G$).
Then $F^\prime(x) = cG(x)$ so that 
\[ \Bigl(F(f(x))\Bigr)^\prime = cG(f(x))f^\prime(x) = cG(x)=F^\prime(x)\]
and so $F(f(x))=F(x)+b$ for some constant $b$. But $F(f(a))=1=F(a)+1$, so $b=1$.
\nl

\end{document}